%% file: thintubes.tex
\documentclass[11pt,a4paper,leqno]{amsart}

\include{packages}

\usepackage[totalwidth=15cm,totalheight=24cm]{geometry}
\usepackage{xcolor}

\newcommand{\cA}{\mathcal A}
\renewcommand{\T}{\mathcal T}

\newcommand{\cML}{\mathcal{ML}}
\renewcommand{\bar}{\overline}

\newcommand{\re}{\mathrm{Re}}
\newcommand{\gr}{\mathrm{gr}}

\newcommand{\cT}{\mathcal T}

\renewcommand{\T}{\mathcal T}

\renewcommand{\bar}{\overline}
\renewcommand{\M}[1]{{\mathcal M}\left(#1\right)}

\newcommand{\argsinh}{{\mathrm{arsinh}}}

\newcounter{notes}%

\begin{document}

\title[]{Filling Riemann surfaces by hyperbolic Schottky manifolds of negative volume} 

\author{Tommaso Cremaschi}
\address{Tommaso Cremaschi:
Trinity College Dublin, School of Mathematics,
17 Westland Row, Trinity College Dublin, Dublin 2, Ireland}
\email{cremasct@tcd.ie}

\author{Viola Giovannini}
\address{Viola Giovannini:
University of Luxembourg, Department of Mathematics, 
Maison du nombre, 6 avenue de la Fonte,
L-4364 Esch-sur-Alzette, Luxembourg}
\email{viola.giovannini@uni.lu}

\author{Jean-Marc Schlenker}
\address{Jean-Marc Schlenker:
University of Luxembourg, Department of Mathematics, 
Maison du nombre, 6 avenue de la Fonte,
L-4364 Esch-sur-Alzette, Luxembourg}
\email{jean-marc.schlenker@uni.lu}

%\thanks{}
%
%\date{v0, \today}

\maketitle

\begin{abstract}
We provide conditions under which a Riemann surface $X$ is the asymptotic boundary of a  convex co-compact hyperbolic manifold, homeomorphic to a handlebody, of negative renormalized volume. We prove that this is the case when there are on $X$ enough closed curves of short enough hyperbolic length.
\end{abstract}

\section{Introduction and results}

\subsection{Hyperbolic manifolds of smallest volume}

The volume of a closed hyperbolic 3-manifold can be considered as a measure of its ``complexity'', and it is natural to ask what is the closed, orientable hyperbolic manifold of smallest volume. The answer is the Weeks manifold \cite{gabai-meyerhoff-milley}.

Consider now a compact Riemann surface $X$. We can extend the previous question in the following manner -- the case of closed hyperbolic manifolds corresponds to $X=\emptyset$.

\begin{question} \label{q:minimal}
 Given $X$, what is the convex co-compact hyperbolic manifold $M$ of smallest volume, with asymptotic boundary $X$?
\end{question}

Convex co-compact hyperbolic manifolds have infinite volume, but they have a well-defined {\em renormalized volume} (see Section \ref{ssc:renormvol} below) which we consider here. The notion of renormalized volume was introduced first in the physics literature by Skenderis and Solodukhin \cite{skenderis-solod}, and then quickly introduced in the mathematics study of conformally compact Einstein manifolds \cite{graham-witten}. For 3-dimensional hyperbolic manifolds, it is closely connected \cite{Holography} to the Liouville functional studied e.g. in \cite{TZ-schottky,Takhtajan:2002cc}. More recently, an explicit upper bound on the renormalized volume of quasifuchsian manifolds in terms of the Weil-Petersson distance between the conformal metrics at infinity, as well as a bound on the difference between the renormalized volume and the volume of the convex core \cite{compare}, led to bounds on the hyperbolic volume of mapping tori \cite{kojima-mcshane,brock-bromberg:inflexibility2}. Moreover, the study of the gradient flow of the renormalized volume has brought a number of new results, see e.g. \cite{bridgeman-bromberg-pallete:convergence,bridgeman-brock-bromberg:gradient,bridgeman-brock-bromberg}. 

Beyond those mathematical motivations, Question \ref{q:minimal} also occurs naturally from a physical perspective, and specifically from the AdS/CFT correspondence. Very briefly, the AdS/CFT correspondence asserts the equality between the partition function of a conformal field theory (CFT) on a $d$-dimensional manifold $X$ and a sum, over all $d+1$-dimensional manifolds $M_i$ with boundary $X$, of a function of the action of a certain (super-)string theory on $M_i$. In a certain ``gravity'' limit, where many features disappear, it reduces to a very special and simplified statement: given a Riemann surface $X$, the partition function of a certain CFT on $X$ should be recovered as a sum of exponential of minus a constant times the renormalized volumes of all convex co-compact hyperbolic manifolds $M_i$ having $X$ as asymptotic boundary:
$$ \cA(X) = a_0 \sum_{\partial M_i=X} e^{-cV_R(M_i)}~. $$
where $a_0$ and $c$ are constants.
In this simplified view, the main term on the $d+1$-dimensional ``bulk'' side corresponds to the convex co-compact manifold $M_i$ with the smallest renormalized volume.

This AdS/CFT correspondence leads to some conjectural statements. For instance, if $X$ is disconnected, the CFT should behave independently on the two connected component, and it might therefore be expected that the convex co-compact manifold of smallest volume ``filling'' $X$ should also be disconnected (see \cite{averages} for a more elaborate analysis).

For instance, if $X=X_+\cup X_-$ is the disjoint union of two connected Riemann surfaces of genus at least $2$, with $X_-$ equal to $X_+$ with opposite orientation, we can compare:
\begin{itemize}
\item the Fuchsian manifold $M_F$ with ideal boundary $X_+\cup X_-$, which has (with the normalization used here) renormalized volume zero,
\item any possible filling of $X_+\cup X_-$ by the disjoint union of two handlebodies $M_+$ and $M_-$, with $\partial_\infty M_+=X_+$ and $\partial_\infty M_-=X_-$.
\end{itemize}
The heuristics above suggests that one of the disconnected fillings might have negative renormalized volume. This might be a motivation for the following conjecture, attributed to Maldacena (see \cite{VP2019}).

\begin{conjecture} \label{cj:maldacena}
  Any connected Riemann surface of genus at least $2$ is the asymptotic boundary of a Schottky manifold of negative renormalized volume. 
\end{conjecture}

By ``Schottky manifold'' here we mean a convex co-compact hyperbolic manifold homeomorphic to a handlebody. 

\subsection{Results}

In what follows $S$ will always denote a closed orientable surface of genus at least $2$. 

\subsubsection{Existence of fillings of minimal renormalized volume}

Before we consider the questions above, it is useful to know that, given a Riemann surface $X$ of finite type, there is at least one convex co-compact filling of $X$ of minimum renormalized volume, and that the set of those minimum volume fillings is finite. Precisely, let $\M X$ be the set of convex co-compact hyperbolic manifolds with ideal conformal boundary $X$, then we think the following question should have a positive answer.

\begin{question}
  Let $V\eqdef \inf_{M\in \M X} V_R(M,X)$. There exists $M_V\in \M X$ such that $V_R(M_V)=V$.
\end{question}

This extends to the case of closed hyperbolic manifold, when $X$ is empty, in which the Weeks manifold is the unique smallest volume closed hyperbolic 3-manifold \cite{Weeks}.

\subsubsection{An upper bound on the renormalized volume}

The main result here is an upper bound on the renormalized volume of a Schottky manifold, when it is obtained from a pants decomposition for which some of the curves are short. We denote by $\epsilon_0$ the 2-dimensional Margulis constant, equal to $\epsilon_0=2\argsinh(1)$. Given a pants decomposition $P$ of $S$, we denote by $M_P$ the handlebody with boundary $S$ in which all curves of $P$ are null-homotopic in $M_P$ (see \S\ref{sssc:handlebodies}), and by $M_P(X)$ the convex co-compact hyperbolic manifold homeomorphic to $M_P$ with complex structure at infinity $X$. The complex structure has a unique hyperbolic metric in its conformal class and we will take lengths with respect to that. Thus, by $\ell_X(\gamma)$ we mean the length of $\gamma$ with respect to the hyperbolic structure induced by $X$. In the case in which there is no ambiguity we will often use $\ell(\gamma)$.

\bthm \label{mainthm}
Let $X$ be a closed Riemann surface of genus $g\geq 2$. Assume that there are $k$ disjoint simple closed curves $\gamma_1, \cdots, \gamma_k$ such that $\ell(\gamma_i)\leq 1, 1\leq i\leq k,$ and there are no other geodesic loops of length less or equal than $1$ in $X$. Then there exists a pants decomposition $P$ containing the $\gamma_i$'s such that
\[ V_R(M_P(X))\leq -\f{\pi^3}{\sqrt e}\sum_{i=1}^k \frac 1{\ell(\gamma_i)} + \left(9+\frac 3 4\coth^2\left(\f14\right)\right)k + 81\coth^2\left(\f14\right)\pi(3g-3-k)(g-1)^{2}~. \] 
\ethm

By imposing the right hand side of the estimate in Theorem \ref{mainthm} to be negative we obtain for instance the following corollary.

\bcor\label{maincoro}
  For all $g\in \N$ s.t. $g\geq 2$, $0<k\leq 3g-3$ and $0<k_1\leq k$ there exists an explicit constant $A=A(g,k_1,k-k_1)>0$ such that if $X$ is a Riemann surface with $k_1$ geodesic loops of length less than $A$ and $k$ geodesic loops of length at most $1$, then $X$ admits a Schottky filling with negative renormalized volume.
\ecor

\brem\label{edgecaseseample}
Let us see a couple of examples for Corollary \ref{maincoro} in the two limit cases.
\begin{itemize}
    \item \textbf{Case $k=k_1=1$.} By the inequality of Theorem \ref{mainthm}, we have    \[ A(g,1)< \dfrac{\pi^3}{\sqrt{e}(9+\f 34 \coth^2(1/4) + 81\coth^2(1/4)\pi(3g-4)(g-1)^2)}~,\] 
    which  is largest for genus $g=2$, in which case the bound is
    \[A(2,1)\leq 0.00221.\] 
    Note that, for large genus $g$, we obtain the asymptotic $A(g,1)\sim cg^{-3}$, for a $c>0$ that can be read from the expression of $A(g,1)$. 

    \item \textbf{Case $k=k_1=3g-3$.} By Theorem \ref{mainthm}, since $k-k_1=0$, we are looking for an $A\eqdef A(g, 3g-3)$ such that
    \[A< \dfrac{\pi^3}{\sqrt{e}(9+\f 3 4\coth^2(1/4))}~,\] 
    we can then take
    \[A(g, 3g-3)= 0.87458~.\]
    
\end{itemize}
\erem

This statement can be compared to \cite[Corollary 5.6]{MD2019}, also see \cite[Theorem 2.1]{VP2019}, which states that: if a Riemann surface $X$ of finite type and genus $g\geq 2$ has  $g-1$ closed curves $\gamma_1, \cdots, \gamma_{g-1}$ such that the complement of their union is a disjoint union of k-holed tori, and if
$$ \frac 1 {\pi-2}\left(\sum_{i=1}^{g-1}\sqrt{\ell_X(\gamma_i)}\right)^2\leq \pi(g-1)~, $$
then
\[ V_R(X,\mathcal P)\leq \pi (g-1)\left(3-\frac{\pi(\pi-2)(g-1)}{\left(\sum_{i=1}^{g-1}\sqrt{\ell_X(\gamma_i)}\right)^2} \right)\]
which is negative if
$$ \sum_{i=1}^{g-1}\sqrt{\ell_X(\gamma_i)}\leq\left(\frac {\pi(\pi-2)(g-1)}3\right)^\frac 12~, $$
and, in the case $g=2$, leads to a better $A(2,1)=\f{\pi(\pi-2)}{3}$.

\subsubsection{Outline of the proof}

The proof of Theorem \ref{mainthm} follows several steps. First, we introduce in Section \ref{sc:3} a notion of ``symmetric'' Riemann surfaces -- those which admit an orientation-reversing involution with quotient a surface with boundary. We prove that given any Riemann surface $X$ of finite-type and any pants decomposition $P$ of $X$, there is a symmetric surface $X_s$ (for which the involution leaves $P$ invariant component-wise) obtained from $X$ by earthquakes along the curves of $P$ (see Lemma \ref{FNtwist}).

Then, in Section \ref{sc:symmetric}, we estimate the renormalized volume of ``symmetric'' Schottky fillings of symmetric surfaces. In Section \ref{sc:variation}, we provide a formula for the difference of the renormalized volume of a filling under an earthquake path of the boundary surface (see Theorem \ref{dehntwistpath}). The result expresses the estimates in terms of the Schwarzian derivative at infinity (see  Section \ref{ssc:schwarzian}) at the core of tubes associated to the pants curves. Finally, Section \ref{sc:proofs} contains the proofs of the main results.

\subsubsection{Convex co-compact fillings}

The result in bounding the difference of renormalized volume under earthquake, see Theorem \ref{dehntwistpath}, can also be applied in the more general setting of convex co-compact manifolds. Specifically it makes sense in the setting where $N(X_0)\in CC(M)$ is a convex co-compact hyperbolic 3-manifold, homeomorphic to $M$, with conformal boundary $X_0\in \mathcal T(\partial M)$. In this more general setting the boundary of $M$ can be disconnected and can be decomposed as $\partial M= F_c\cup F_i$ where $F_i$ does not compress in $M$ and each component of $F_c$ compresses (i.e. it has at least a loop bounding a disk in $M$).

 Let $c_{\mathbf t}^{\mathfrak m}:[0,1]\rar CC(M)$ be an earthquake path (we earthquake by a parameter $t_i$, with $\mathbf t = (t_1, \dots, t_n)$, along the curve $\gamma_i)$ along a multi-curve $\mathfrak m=\set{\gamma_i}^n_{i=1}\subset S$ which, with respect to the reference hyperbolic metric $X_0$, can be subdivided into:
\begin{itemize}

\item $\mathfrak m_{1}^c$: the set of { compressible geodesic loops $\gamma$ of $\mathfrak m$ with length at most $1$};
\item $\mathfrak m_{ 1}$:  { the set of geodesic loops $\gamma$ in $F_c$ and not in $\mathfrak m_{1}^c$ such that any compressible geodesic loop $\alpha$ intersecting $\gamma$ essentially has length at least $1$;}
\item $\mathfrak m_\infty$: the set of geodesic loops $\gamma$ of $\mathfrak m$ that are { contained in $F_i$ and so incompressible}.

\end{itemize}
 { Note that not every $\m$ admits such a decomposition with respect to the given $X_0$, as there could be a $\gamma_i\in \mathfrak{m}$ in a compressible component, of length more than $1$ and intersecting a short compressible loop. }

\bthm\label{dehntwistpath2} Let $X_0\in\T(\partial M)$ and $\m=\m_\infty\cup \m_{1} ^c\cup \mathfrak m_{1}$ be a multi-curve and $c_{\mathbf t}^{\mathfrak m} $ be an earthquake path terminating at $X_1$. Then \small
\[ \abs{ V_R(X_1) -V_R(X_0) }\leq   \sum_{\gamma_i\in \pi_0(\mathfrak m_{1}^c)} (3\ell_i\coth^2 \left( \ell_i/4\right) ) t_i +C \sum_{\alpha_j\in\mathfrak m_{ 1} } t_j\ell_j +{3} \sum_{\beta_k\in\mathfrak m_\infty}  t_k \ell_k~,\] 
for  $C=3\coth^2\left(\frac {1} 4\right)<50.013$.\normalsize

\ethm

\subsection*{Acknowledgements}
TC was partially supported by MSCA grant 101107744--DefHyp. VG was partially supported by FNR AFR grant 15719177. JMS was partially supported by FNR OPEN grant O20/14766753. We also thank the referee for their mindful comments that helped increase the readability of the paper.

\section{Notation and background}
In this section we recall the main objects and tools that we will use in this work.

\subsection{Hyperbolic surfaces}\label{hypmanifolds}

\subsubsection{Teichm\"uller space} Good references for Teichm\"uller space are \cite[Chapter 6-7]{Hubbard2016} and \cite[]{FM2011}, we now recall what we will need.
Any closed, oriented, surface of genus $\geq 2$ is hyperbolic, i.e. is homeomorphic to a quotient $\mathbb{H}^2/\Gamma$ of the hyperbolic space by a discrete, torsion-free subgroup of the orientation-preserving isometries of $\mathbb{H}^2$. The Teichm\"uller space of $S$ can be defined in the following various ways, {depending on the set-up, we will use the most suitable definition}: 
\begin{align*}
\mathcal{T}(S) & =\{\text{h hyperbolic metric on }S\}/\text{Diffeo}_0(S)~, \\
\mathcal{T}(S) & =\{c\text{ complex structure on }S \}/\text{Diffeo}_0(S)~, \\
\mathcal{T}(S) & =\{[g] \text{ s.t. g is a Riemannian metric on S}\}/\text{Diffeo}_0(S)~. \end{align*}
Here $\text{Diffeo}_0(S)$ is the group of diffeomorphisms of $S$ isotopic to the identity, and it acts by pull-back, moreover $g_1\in [g_2]$ if and only if there exists a {smooth function $u_1\colon S\rightarrow \mathbb{R}$} such that $g_2=e^{u_1}g_1$, i.e. $[g]$ represents the class of Riemannian metrics conformal to $g$. In particular, to any complex structure on $S$ corresponds a conformal class of metrics $[g]$, in which, by the Riemann uniformization Theorem \cite[Theorem 1.1.1]{Hubbard2016}, there exists a unique hyperbolic metric $h\in[g]$.

\subsubsection{Margulis tubes}\label{ssc:margtubes}

Hyperbolic surfaces, and in general hyperbolic $n$-manifolds, have the important property that ``short geodesics" have particularly nice neighbourhoods.  We will mostly deal with the surface case and so we restrict ourselves to that setting.

\bdefi \label{margulitubedefi}
By a \emph{thin tube}, for a hyperbolic surface, we mean the set of points $\mathbb T(\ell)$ around a geodesic $\gamma$ of length $\ell\leq \epsilon_0$, { with $\epsilon_0=2\argsinh(1)$ the $2$-dimensional Margulis constant}, that are at a distance at most $L\eqdef\argsinh\left(\f1{\sinh\left(\f\ell 2\right)}\right)$. In $\mathbb T(\ell)$ the injectivity radius is bounded as 

\[\f{\ell}2 \leq \operatorname{inj}(p)=\argsinh\left(\sinh(\ell/2)\cosh(L-d)\right),\qquad d=d(p,\partial \mathbb T)~,\] 
and its maximum is achieved on $\partial \mathbb T(\ell)$, see \cite[Thm 4.1.6]{Bu1992}. The hyperbolic metric on $\mathbb T(\ell)$ can be written as $d\rho^2+\left(\f{\ell}{2\pi}\right)^2\cosh^2(\rho)d\theta^2$, $\theta\in[0,2\pi]$ and $\rho\in[-L,L]$. Moreover, any multi-curve $P$ such that each component is simple and has length at most $\epsilon_0$ can be completed to a pants decomposition of $S$. For details see \cite[Thm 4.1.1]{Bu1992}.
\edefi

\subsubsection{Earthquakes along simple closed geodesics} 

We recall here some basic facts on earthquakes along closed geodesics, which will be needed. For more background see \cite[Sec 10.7.3]{FM2011} and \cite[Part III]{CEM2006}.

Given a simple closed geodesic $\gamma$ on a hyperbolic surface $(S, h)$ a (left) $t$-earthquake is a map $\phi_{\gamma,t}$ from $S$ to itself, discontinuous along $\gamma$, defined by cutting $S$ along $\gamma$, twisting the left-hand side of $\gamma$ by a fixed length $t$ in the positive direction, and gluing back isometrically the two sides.\footnote{Note that this definition requires the choice of an orientation of $\gamma$, but the result does not depend on which orientation is chosen.}

Taking the push-forward of the hyperbolic metric by $\phi_{\gamma,t}$ defines a new hyperbolic metric on $S$, and in this manner $\gamma$ and $t$ define a homeomorphism of $\cT(S)$, which is also called the left earthquake of length $t$ along $\gamma$, and denoted by {$E_{\gamma}(t)$}. 

By continuously varying the twisting length $t$, one gets a path of diffeomorphisms of $S\setminus\gamma$, and by pulling back $h$ through such a path, we get a path in $\mathcal{T}(S)$. 

 Let us now define earthquakes more carefully. Having fixed a simple closed curve $\gamma$ in $S$, we consider the unique geodesic on $(S,h)$ in the same isotopy class again by $\gamma$. {In this way, the operation only depends on the isotopy class of $\gamma$.} Let $\ell$ be the length of $\gamma$ with respect to $h$, and $N_r\cong S^1\times [-r,r]\cong\mathbb{R}/\ell\mathbb{Z}\times [-r,r]$ be the tubular $r$-neighborhood of $\gamma$ parameterized in such a way that $S^1\times \{0\}$  isometrically identifies with $\gamma$ and $\{e^{i\theta}\}\times [-r,r]\in S^1\times [-r,r]$ with a { geodesic} segment of length $2r$ orthogonal to $\gamma$ parameterized in unit velocity by the coordinate in $[-r,r]$. We now choose an arbitrary function $f\colon [-r,r]\rightarrow \mathbb{R}$ such that $f$ is smooth on $[-r,r]\setminus\{0\}$, increasing on $[-r,0]$, constantly equal to $0$ in a neighborhood of $-r$, constantly equal to $1$ in a left neighbourhood of $0$, and equal to $0$ on $(0,r]$. We then define $\phi_{\gamma,t}\colon S\rightarrow S$, with $t\in[0,\infty)$, as the diffeomorphism of $S\setminus\gamma$ such that $\phi_{\gamma,t}$ is the identity outside $N_r$, and 
\[\phi_{\gamma,t}(e^{i\theta}, r)=\left(e^{i(\theta+\frac{2\pi}{\ell}tf(r))},r\right)\] 
for any $(e^{i\theta}, r)\in N_r$. Note that $\phi_{\gamma,0}$ is the identity, and that $\phi_{\gamma,\ell}$ extends to a diffeomorphism of $S$, which is called a \textit{Dehn twist}.

As $\phi_{\gamma,t}$ acts by isometry on $N_\epsilon(\gamma)\setminus \gamma$ and  fixes the metric on $\gamma$, the push-forward $(\phi_{\gamma,t})_{*}(h)$ is a new well defined hyperbolic Riemannian metric on $S$. We say that $(\phi_{\gamma,t})_{*}(h)$ is obtained by a \textit{(left) earthquake of parameter t} along $\gamma$. 

We define $\phi_\gamma:[0,a]\rar \mathcal T(S)$, $a>0$, to be a \emph{earthquake path along $\gamma$} by $\phi_\gamma(t)=(\phi_{\gamma,t})_{*}(h)$. The \emph{infinitesimal earthquake along $\gamma$} is the derivative of $\phi_\gamma$ in $t$ at $t=0$, this can also be seen as a vector field $v$ on $S$ by differentiating the path of diffeomorphisms ${(\phi_{\gamma, t})}_{t\in [0,\epsilon]}$  with respect to $t$ and evaluating it at $t=0$. For more background see \cite[Sec 10.7.3]{FM2011} and \cite[Part III]{CEM2006}.

\subsection{Hyperbolic 3-manifolds.}

Some references on hyperbolic 3-manifolds are \cite{Mar16,MT1998,He1976}, we now recall what we will need in this work. A 3-manifold $M$ is hyperbolic if it is homeomorphic to $\mathbb H^3/\Gamma$ for $\Gamma$ a discrete, torsion free subgroup of $\mathbb{P}SL(2, \mathbb{C})$, the positive isometry group of $\mathbb H^3$. 

The action of $\Gamma$ on $\mathbb{H}^3$ can be naturally extended to $\partial\mathbb{H}^3=\hat{\mathbb{C}}$, with $\hat{\mathbb{C}}$ the Riemann sphere, but it does not remain properly discontinuous, that is, the closure of the orbit of a point $x\in \mathbb{H}^3$ has non-empty set of accumulation points $\Lambda_x(\Gamma)$ in $\mathbb{H}^3\cup \partial\mathbb{H}^3$. One can show that actually $\Lambda_x(\Gamma)$ does not depend on $x$. We then denote it simply by $\Lambda(\Gamma)$ and we call the complement $\Omega(\Gamma)=\partial\mathbb{H}^3\setminus\Lambda(\Gamma)=\hat{\mathbb{C}}\setminus\Lambda(\Gamma)$ the \textit{domain of discontinuity of $\Gamma$}. We observe that $\Lambda(\Gamma)$ is closed, and that both $\Lambda(\Gamma)$ and $\Omega(\Gamma)$ are $\Gamma$-invariant. The action of $\Gamma$ on $\Omega(\Gamma)$ is properly discontinuous, we can then define the \textit{boundary at infinity of $M=\mathbb{H}^3/\Gamma$} as the surface 
\[\partial_{\infty}M=\Omega(\Gamma)/\Gamma~.\] 
Since $\Omega(\Gamma)$ is an open subset of $\hat{\C}$, and the elements of $\mathbb{P}SL(2,\mathbb{C})$ are in particular bi-holomorphism of ${\hat{\mathbb{C}}}$, the boundary at infinity $\partial_{\infty}M$ of $M$ is naturally equipped with {a \textit{complex projective structure}, and thus also a complex structure.}
{A complex projective structure on a surface $S$ is an atlas of charts to $\hat{\mathbb{C}}$ whose transition maps are restriction of M\"obius transformations. Equivalently, a complex projective structure is the datum of a holonomy representation of $\pi_1(S)$ in $\mathbb{P}SL(2,\mathbb{C})$, and an equivariant \textit{developing map}, that is, an immersion of the universal cover $\widetilde{S}$ equipped with the lifted complex projective structure in $\hat{\mathbb{C}}$, which locally restricts to projective charts. The developing map is unique up to composition with M\"obius transformations. The deformation space of complex projective structures forms a holomorphic vector boundle $\pi$ on the Teichm\"uller space of $S$, of dimension $12g-12$ (see \cite{dumas-survey}). Given $X\in \mathcal{T}(S)$, the fiber $\pi^{-1}(X)$ is parameterized by the Schwarzian derivative of the developing map of each point in the fiber (see Section \ref{ssc:schwarzian}).}

\subsubsection{The convex core}\label{ssc:concore} We define the \textit{convex core} of $M=\mathbb{H}^3/\Gamma$ as 
\[ C(M)=\text{Hull}(\Lambda(\Gamma))/\Gamma~, \] 
where $\text{Hull}(\Lambda(\Gamma))$ is the convex envelop of the points of $\Lambda(\Gamma)$ in $\mathbb{H}^3\cup \partial\mathbb{H}^3$.
The convex core of $M$ is also characterized as the smallest non-empty \textit{strongly geodesically convex}\footnote{Here we say that a subset $K\subset M$ is ``strongly geodesically convex'' if any geodesic segment in $M$ with endpoints in $K$ is entirely contained in $K$.} subset of $M$, that is, the smallest convex subset of $M$ which is also homotopically equivalent to $M$. It is also not difficult to prove that if $M$ has finite volume, then the limit set $\Lambda(\Gamma)$ coincides with { $\hat \C$}, and so $C(M)=M$. Here, we will be interested in the case of $M$ having infinite volume. The convex core $C(M)$ is generically a 3-dimensional domain, but in some cases, it can be a totally geodesic surface in $M$, possibly with geodesic boundary.

\bdefi
A hyperbolic $3$-manifold $M=\mathbb{H}^3/\Gamma$ is \textit{convex co-compact} if its convex core $C(M)$ is compact.
\edefi

When $M$ is convex co-compact then $\overline M=M\cup (\partial_{\infty}M)$ is its manifold compactification and $\partial_{\infty}M$ is homeomorphic { to} the closed surface $S=\partial\bar{M}$, and so \[[\partial_{\infty}M]\in \mathcal{T}(S)~.\]

We call \textit{end} a connected component of $M\setminus C(M)$, or, more generally, of the complement of a strongly geodesically convex compact subset of $M$. An end is homeomorphic to $S^i\times [0, +\infty)$, with $S^i$ a connected component of the boundary $S=\partial_{\infty}M$, and it has infinite hyperbolic volume.
  
We denote by $CC(M)$ the space of convex co-compact hyperbolic structures on $M$ considered up to homotopy equivalence. The deformation space $CC(M)$ is parameterized by the one of conformal structures on the boundary at infinity, see \cite[Thm 5.1.3.]{Ma2016} and \cite[Thm 5.27]{MT1998}:
\[ CC(M)=\quotient{\mathcal T(\partial\overline M)}{T_0(D)},\]
where $T_0(D)\subset MCG(\partial\overline M)$ is the subgroup generated by Dehn twists along compressible curves\footnote{{An essential loop $\gamma$ in $\partial \overline M$ is compressible if it is null-homotopic in $M$, i.e. it bounds a compressing disk in $M$.}} of $\partial\overline M$ and $\mathcal T(\partial\overline M)$ is the product of the Teichm\"uller spaces of the connected components of $\partial \bar{M}$.

\subsubsection{Handlebodies.} \label{sssc:handlebodies}

We will think of an handlebody $H_g$ of genus $g\geq 1$ as the following data. Given a surface $S=S_g$ and a pants decomposition $P$ on $S$ we can form the $3$-manifold $H_0$ by attaching $3g-3$ thickened disks $\mathbb D^2\times I$ to $S\times I$ by gluing each $\partial \mathbb D^2\times I$ to $N_\epsilon (\gamma) \times\set 0$ for $\gamma\in P$. The manifold $H_0$ has then a genus $g$ boundary component and ${2g-2}$ sphere boundary components. After filling each sphere component with a 3-ball we obtain a handlebody $H_P\cong H_g$, this is unique up to isotopy. We will think of this as the handlebody induced by $P$. Note that $H_P$ is well defined up to isotopy.

 We now define what it means to \emph{fill} a given conformal structure $X\in\mathcal T(S)$ via a (complete) hyperbolic 3-manifold $M$ so that $M$ is homeomorphic to a handlebody and its conformal boundary is $X$.
 
\bdefi\label{fillingdefi} Given a conformal structure $X\in\mathcal T(S_g)$ and a pants decomposition $P$ on $S_g$ we say that $M_P(X)$ is the \emph{Schottky filling} of $X$ with pants curve $P$ if it is the (complete) hyperbolic 3-manifold obtained by uniformising ${H_P}$ so that its conformal boundary is $X$. By $CC_P(S_g)$ we denote the deformation space of a hyperbolic genus $g$ handlebody obtained by gluing disks along $P$.
\edefi

\brem\label{generalhandlebody}
More generally, a handlebody $H_g$, topologically, is any irreducible compact 3-manifold $M$ with a unique boundary component such that {the map induced by the inclusion $\partial M \hookrightarrow M$ on the fundamental groups is surjective}, \cite[]{He1976}. Thus, the manifold $M\eqdef F\times I$, for $F$ a compact orientable surface with non-empty boundary, is also a handlebody with boundary given by the double of $F$ along $\partial F$. In the case that $F$ is not-orientable then we can consider the twisted $I$-bundle\footnote{{Recall that a twisted $I$-bundle is a non-trivial $I$-bundle, i.e. $N\not\cong F \times I$ .}} $N=F\overset{ \sim}{ \times }I$ in which $\partial N$ is given by the orientation double cover of $F$.
\erem

\subsubsection{The Schwarzian derivative at infinity} \label{ssc:schwarzian}

Given a Riemann surface, a holomorphic quadratic differential is a holomorphic section of the symmetric square of its holomorphic cotangent bundle, and in holomorphic coordinate it can be expressed as $\varphi(z)dz\otimes dz=\varphi(z)dz^2$.
Let $D\subseteq \mathbb{C}$ be a connected open set, and $f\colon D\rightarrow \hat{\mathbb{C}}$ a locally injective holomorphic map. The \textit{Schwarzian derivative} of $f$ is the holomorphic quadratic differential \[\mathcal S(f)=\biggl(\biggl(\dfrac{f^{''}}{f^{'}}\biggr)^{'}-\dfrac{1}{2}\biggl(\dfrac{f^{''}}{f^{'}}\biggr)^{2}\biggr)dz^2~.\]
The Schwarzian derivative has the following properties:
    \begin{enumerate}
        \item Let $f$ and $g$ be locally injective holomorphic maps such that the composition is well defined, then \[\mathcal S(f\circ g)=g^{*}\mathcal S(f)+\mathcal S(g)~.\]
        %% \item For any $f\in \mathbb{P}SL(2, \mathbb{C})$ \[S(f)=0~.\]
        \item For any holomorphic map $f:U\to \C$, where $U\subset \C$ is an open subset, $\mathcal S(f)=0$ if and only if $f\in \mathbb{P}SL(2, \mathbb{C})$, that is, if and only if $f$ is the restriction to $U$ of a M\"obius transformation.
    \end{enumerate}

We will be interested in considering the Schwarzian derivative of {the developing map $f\colon \mathbb{H}^2\rightarrow \Omega(\Gamma)$ of a non simply connected domain of discontinuity $\Omega(\Gamma)$, whose quotient by $\Gamma$ gives the boundary at infinity $\partial_{\infty}M=\mathbb{H}^2/\Gamma'$ of $M=\mathbb{H}^3/\Gamma$. The map $f$ is a covering of $\Omega({\Gamma})$, hence locally univalent, and it is $(\Gamma',\Gamma)$-equivariant, thus, thanks also to property $(2)$, the Schwarzian $\mathcal S(f)$ of $f$ descends to a holomorphic quadratic differential on $\partial_{\infty}M$.}

\subsubsection{The renormalized volume} \label{ssc:renormvol}

If one is willing to talk about volumes for convex co-compact hyperbolic structures on $M$, being this infinite, some kind of renormalization will be needed. A possibility is to consider the function 
\[V_{C}\colon CC(M)\longrightarrow \mathbb{R}_{\geq 0}~,\]
which associates to any convex co-compact structure $M$ the volume of its convex core $\text{Vol}(C(M))$. %It was was proven by Bridgman and Canary that $V_{C}$ is continuous \cite{BC2017}. 
The renormalized volume is some kind of relative of the function $V_{C}$, which presents much better analytic properties.

The idea is to consider an exhaustion of $M$ by strongly geodesically convex compact subsets $\{C_r\}_r$ coming together with an equidistant foliation of the ends, and to renormalize the associated volumes $\text{Vol}(C_r)$ in order to get a finite number which does not depend on $r$.

Before giving the definition of renormalized volume, we need to introduce some preliminary notions. 

\bdefi
    Let $M$ be convex co-compact and let $C\subseteq M$ be a compact, convex subset with smooth boundary. We define the {\em $W$-Volume of $C$} as 
    \[W(C)=\text{Vol}(C)-{\dfrac{1}{2}}\int_{\partial C} H dA_{\partial C}~,\]
    where $\text{Vol}(C)$ is the hyperbolic volume of $C$ with respect to the metric of $M$, $H$ is the mean curvature of $\partial C$, and $dA_{\partial C}$ is the area form of the induced metric on the boundary $\partial C$.
    \edefi
    
The \textit{mean curvature} is {half} the trace of the \textit{shape operator} $B(X)=\nabla_{X}(N)$ with $N$ the { outer} unit normal to $\partial C$ and $\nabla$ the Levi-Civita connection, and any vector field $X\in T(\partial C)$.

{In what follows, we assume the compact subset $C$ to be strongly geodesically convex, so that it is homotopically equivalent to $M$ and can be used to decompose the manifold $M$ in neighborhoods of its convex co-compact ends and a compact piece containing the convex core $C(M)$}. The additional term with the mean curvature { in the definition just above} is the right one to get a good renormalization. In \cite{S2008}, it is indeed proven that, denoting by $C_r$ the $r$-neighborhood of $C$ in $M$, for any $r\geq 0$ 
 \begin{equation}\label{Wvol} W(C_r)+r\pi\chi(\partial_{\infty}M)= W(C) \end{equation}
where $\chi(\cdot)$ is the Euler characteristic.

\bdefi
    Let $E$ be an end of $M$. An \textit{equidistant foliation} is a foliation $\{S_r\}_{r\geq r_0}$ of a neighborhood of $\partial_{\infty}M$ in $\bar{M}$ in convex surfaces, such that for any $r'>r>r_0\geq 0$ the surface $S_{r'}$ lives between $S_r$ and $\partial_{\infty}M$, and its points stay at constant distance $r' -r$ from $S_r$.
\edefi

By definition, the boundaries of the $C_r$'s form an equidistant foliation $\{\partial C_r\}_{r}$ of the ends in $M\setminus C$.

Given any $C\subseteq M$ as above, and any end $E_i=S^i\times [0, +\infty)$ in $M\setminus C$, we can consider the associated equidistant foliation $\{\partial^i C_r\}_r$, for $r\geq 0$, where with $\partial^i$ we mean the connected component of $\partial C_r$ facing $S^i$ in the boundary at infinity $S=\partial_{\infty}M$. Let us call $g_r$ the induced metric on $\partial C_r$. Then we can define a metric on the boundary at infinity as
\begin{equation}\label{I}g:={\lim_{r\rightarrow +\infty}4e^{-2r}g_{r}} \in [\partial_{\infty}M]~,\end{equation}
see \cite[Def. 3.2]{compare}, or \cite[Def. 3.2]{volumes} for a slightly different point of view. A key property of this metric $g$ is that it is in the conformal class at infinity of $M$, that is, it is compatible with the complex structure at infinity of $M$ {(note that this remains true if we change the factor $4$ in (\ref{I})). Vice-versa, up to scaling by a big enough positive constant, any representative in $[\partial_{\infty}M]$ can be realized in this way, \cite{EpsteinSurf}}. This leads to the following bijective correspondence: 

$$ \begin{matrix}
\begin{Bmatrix}
 \text{Riemannian metrics g on } S\\
 \text{such that g}\in [\partial_{\infty}M]\\
  \text{up to multiplication by }s\in\mathbb{R^{+}}\\

\end{Bmatrix}&&\longleftrightarrow&&\begin{Bmatrix}
 \text{Equidistant convex foliations  } \\
 \text{ of a neighborhood of } \partial_{\infty}M\\
 \text{up to } \sim_{\mathcal{F}}
\end{Bmatrix}&& \\
\end{matrix} $$ 
where two such foliations are $\sim_{\mathcal{F}}$-equivalent if and only if they are { equal outside a compact set}; {rescaling by a positive constant a Riemannian metric in $[\partial_{\infty}M]$ corresponds to a reindexing of the associated foliation}. 

\bdefi
    We define the $W$-Volume of $M$ with respect to $g\in[\partial_{\infty}M]$ as 
    \[W(M,g)=W(C_r(g))+\pi r \chi(\partial_{\infty}M)~,\]
     where $\{C_r(g)\}_{r\geq r_0}$, with $r_0$ big enough, is the exhaustion in compact strongly geodesically convex subsets defined by the equidistant foliation associated to $g$, {indexed in such a way that the sequence of induced metrics $g_r$ on $\partial C_r(g)$ satisfies (\ref{I})}.
\edefi

 Thanks to equation \eqref{Wvol} above, the $W$-volume $W(M,g)$ is well defined.
 {We also remark that
choosing the factor $4$ for the corresponding metric at infinity as in (\ref{I}), is the scaling which makes the geodesically convex subset $C_0$ associated to the hyperbolic metric in the conformal boundary at infinity of a Fuchsian manifold to be isometric to the induced metric on its $2$-dimensional convex core.}

We can finally define the renormalized volume of $M$.

\bdefi
    Given a convex co-compact hyperbolic $3$-manifold $M\in CC(N)$, its {\em renormalized volume} is defined as 
    \[V_R(M)=W(M,h)~,\] 
    with $h\in [\partial_{\infty}M]$ the hyperbolic representative.
\edefi

Thanks to the parametrization of the space of convex co-compact structures $CC(M)$, we can think about the renormalized volume as a function from the Teichm\"uller space:
 \[V_R\colon\mathcal{T}(\partial \bar{M})\longrightarrow \mathbb{R}~.\]

 \brem
 It is possible to define the $W$-volume also for the convex core $C(M)$ of $M$. In this case the boundary $\partial C(M)$ is not smooth, and the integral mean curvature of the boundary is replaced by the length of the {\em measured pleating lamination} (see see \cite{thurston-notes,epstein-marden}):
 \[W(C(M))=\text{Vol}(C(M))-\dfrac{1}{4}L(\beta_{M})~. \]

\erem

The renormalized volume satisfies the following differential formula, see \cite[Corollary 3.11]{compare}. 

\bthm\label{dVRdiff} 
Let $M$ be a convex co-compact hyperbolic $3$-manifold, $\phi_M$ the holomorphic quadratic differential given by the Schwarzian derivative of the developing map of the projective structure of $\partial_{\infty}M$, and $\mu\in T_{[\partial_{\infty}M]}\mathcal{T}(\partial\bar{M})$. Then, the differential of the renormalized volume at $[\partial_{\infty}M]$ satisfies
\[\diff V_R(\mu) =\text{Re} (\langle \mu, \phi_M\rangle).\] 
Here the space of holomorphic quadratic differentials on $\partial_{\infty}M$ is identified with the cotangent bundle $T^{*}_{[\partial_{\infty
}M]}\mathcal T(\partial\bar{M})$ through the Bers embedding, and the pairing $\langle\  ,\  \rangle$ is the duality one with $T_{[\partial_\infty M]}\mathcal T(\partial\overline M)$, which is the space of harmonic Beltrami differentials \cite{Hubbard2016}, \cite{QTT}.
\ethm

\section{Earthquakes to symmetric Surfaces}\label{sc:3}

In this section we study conformal structures { $X$ on a surface $S$ that admit an orientation-reversing involution $\sigma: X\rar X$ such that, if $X$ is equipped with its unique compatible hyperbolic metric, $X_\sigma\eqdef \quotient{X}{\sigma}$ is a hyperbolic surface with totally geodesic boundary}. The main result of this section is Lemma \ref{FNtwist}, which states that given $X\in \T(S)$ and $P\subset S$ a pants decomposition {there exists a symmetric conformal structure $X'\in\mathcal T(S)$ and a path in $CC_P(S)$ from $M_P(X)$ to $M_P(X')$ which is obtained by doing earthquakes of bounded length along the curves of $P$.}

\bdefi
  Let $X\in \T(S)$, then $X$ is a \emph{symmetric surface} if {$S$ admits an orientation reversing involution $\sigma:S\rar S$ that is a local isometry for (the hyperbolic metric on) X} and such that $X_\sigma\eqdef \quotient{X}{\sigma}$ is a surface with non-empty boundary. The subset of Teichm\"uller space of surfaces { for which $\sigma$ is a local isometry} will be denoted by $\T_\sigma(S)$ and the subspace of surface admitting an involution $\sigma$ by $\T_s(S)=\cup_\sigma \T_\sigma(S)$. 
\edefi

\brem
The surface $X_\sigma$ does not have to be orientable.
\erem

\blem\label{symmquotient} Let $X$ be a hyperbolic surface with an orientation reversing involution $\sigma:X\rar X$ that is a local isometry. Then, $\partial X_\sigma=\text{Fix}(\sigma)$ is given by a multi-curve $\mathfrak m\subset X$ such that for each $\gamma\in\pi_0(\mathfrak m)$ we have $\sigma\vert_{\gamma}=\id_\gamma$.
\elem

\bpf
By \cite[Theorem 1.10.15]{klingenberg} the set of fixed points is a closed totally geodesic sub-manifold, thus it is the union of a closed multi-curve $\mathfrak m$ and possibly a finite collection of points. By looking at the action on a small enough ball around an isolated fixed point (so that the centre is the unique fixed point) one can see that, being $\sigma$ orientation reversing, isolated fixed points are not possible and so the fixed set has to be a geodesic multi-curve.

We now want to show that $\mathfrak m$ is the boundary of $X_\sigma$. Let $B\subset X$ be a small enough ball such that $\mathfrak m\cap B$ separates $B$ in two balls and $B=\sigma(B)$. Then, $B/\sigma$ is homeomorphic to a half disk with boundary in $\mathfrak m$. By connectedness and continuity this shows that $\mathfrak m\subset\partial X_\sigma$. The reverse containment follows from the fact that $\sigma:X\setminus\mathfrak m\rar X\setminus \mathfrak m$ is a $2$ to $1$ cover and so $(X\setminus \mathfrak m)/\sigma $ is a surface without boundary. \epf

\begin{figure}[htb!]
\begin{center}
\includegraphics[width=\textwidth]{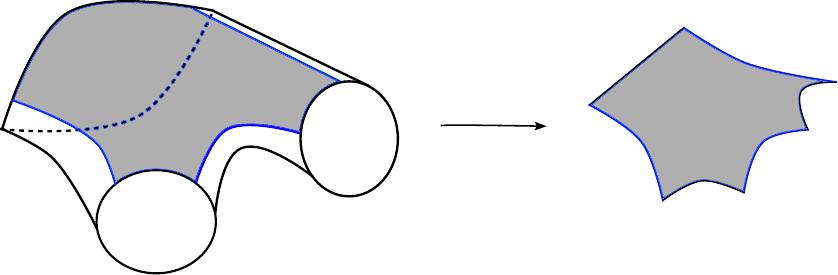}
\caption{The seams (in blue) in a pair of pants with the two hexagons $H_1$ (shaded), $H_2$ and the $\pi_Q$ map.}\label{pantssym}
\end{center}
\end{figure}

In each pair of pants, a seam is the orthogeodesic connecting two distinct boundary components, so each pair of pants has 3 such arcs, see Figure \ref{pantssym}. For every pair of pants $Q$ we have on each boundary component $\gamma_i$ two marked points $x_i^1, x_i^2$, endpoints of the seams of $Q$. We define a marked pants decomposition $P^m$ to be $P$ together with a choice of either $x_i^1$ or $x_i^2$ for each pair of pants $Q$ and each boundary curve of $P$. 

Let $X\in\T (S)$ be a hyperbolic surface, $P$ be a pants decomposition of $X$, and $\mathcal S$ be the set of the induced seams with marked endpoints, i.e. a marked pants decomposition $P^m$. Then, we  define the Fenchel-Nielsen coordinates for $X$ as follows: $FN(X)=(\ell_i,t_i)_{i=1}^{3g-3}$ where the $\ell_i$ are the hyperbolic lengths of the pants curve in the hyperbolic structure on $X$ and the $t_i$ are the twist parameters with respect to the two marked points on the curve $\gamma_i$. The twist parameters are computed by fixing lifts in the universal cover of $\gamma_i$ and then taking their signed euclidean distance. For details on Fenchel-Nielsen coordinates see \cite[Sec 6.2]{Bu1992} or \cite[Sec 10.6]{FM2011}.

 Thus, if $t_i=0$, the seams match up and the two marked points are identified. If $t_i=\ell_i/2$, the seams match up but the marked points are opposite to each other.

Moreover, the seams cut each pair of pants $Q$ into two isometric right-angled hexagons $H_1$ and $H_2$. We can then define an orientation-reversing involution $\sigma_Q:Q\to Q$ which maps $H_1$ to $H_2$ and $H_2$ to $H_1$ isometrically, and is the identity on the seams, see Figure \ref{pantssym}. The quotient of $Q$ by $\sigma_Q$ is then a right-angled hexagon $E_Q$, on which $Q$ projects by a map $\pi_Q:Q\to E_Q$ which is a local isometry outside of the seams.  

\brem\label{seamsglue}
The maps $\set{\pi_Q}_{Q\in P}$ glue together to a map $\pi: X\rar X$ that is an orientation reversing local isometry (outside of the seams) if all seams match up. Moreover, if that is the case then $X_\pi$ is a surface, not necessarily orientable, whose boundary is given, by Lemma \ref{symmquotient}, by the union of the seams.
\erem

\blem\label{FN}
Let $X\in \T(S)$, and let $P=\set{c_1,\cdots, c_{3g-3}}$ be a marked pants decomposition of $X$ and let $(\ell_i, t_i)$ be the corresponding Fenchel-Nielsen coordinates. Then, the Riemann surface $X_0$ with Fenchel-Nielsen coordinates $(\ell_i, t_i')$, $t_i'=0,\ell_i/2$, admits an orientation-reversing isometry which leaves invariant each curve of $P$. 
\elem

\bpf
We want to show that the surface $X_0$ defined by $(\ell_i,t_i')$ admits an orientation-reversing isometry mapping each geodesic loop in $P$ to itself. 
  
The surface $X_0$ is obtained by gluing $3g-3$ pairs of pants with boundary lengths given by the $\ell_i$'s and in the pattern given by $P$ such that if two pairs of pants $Q_1$ and $Q_2$ ($Q_1$ could be equal to $Q_2$) are glued along a geodesic loop $c_i\in \pi_0(P)$ then the endpoint of the seam $y_1\in Q_1\cap c_i$ is glued to $y_2\in Q_2\cap c_i$ without any twist. The pairs of pants $Q_i$, $i=1,2$, are obtained by doubling regular hexagons $E_i$ along the seams, and each $P_i$ is equipped with an orientation-reversing isometry map $\pi_i:Q_i\rightarrow Q_i$ exchanging the two hexagons. The fixed point set of this map is exactly the seams of $Q_i$.

Since the seams on $c_i\subset Q_1\cap Q_2$ have endpoints that are $\ell_i/2$ apart and by our glueing condition one of them matches up we know that they both do. Therefore, all the seams with endpoints on $c_i$ match-up and we can glue the maps $\pi_1$ and $\pi_2$ to obtain an orientation-reversing isometry from $Q_1\cup Q_2$ to $Q_1\cup Q_2$. By doing this for all pants we obtain the required statement.\epf

\brem Given $X$ and $\sigma:X\rar X$ then, for specific markings in the $FN$-coordinates the quotient surface is orientable and equal to a thickening of the glueing graph of the pants decomposition in which if a curve $c_i$ has twist parameter equal to $\ell_i/2$ then the quotient edge is glued with an half-twist.
\erem

\blem\label{FNtwist}
Given a pants decomposition $P$ on $S$ and $X\in\mathcal T(S)$, there exists $X',X_s\in\mathcal T(S)$ such that $X_s$ is symmetric, $M_P(X')\overset{isom}\cong M_P(X)$, and $X_s$ is obtained from $X'$, in $FN_P$ coordinates, by twisting at most $\ell(c_i)/4$ (in the positive or negative direction) over each curve in $P$.
\elem

\bpf First note that in $CC_P(S)$ we can do full twists along curves of $P$ and get isometric structures, see \cite[Thm 5.1.3.]{Ma2016}. Recall that we denote by $M_P(X)\in CC_P(S)$ the structure corresponding to $X\in\mathcal T(S)$ with compressible curves given by $P$.

We use $P$ to define the Fenchel-Nielsen coordinates by choosing seams $y\in\{x_1^i,x_2^i\}\subset c_i$, see Lemma \ref{FN}. Also note that a full twist along $c_i$ has length $\ell_i$. Let $X$ be the given structure, then $FN_P(X)=(\ell_i(X),t_i(X))_{i=1}^{3g-3}$. By doing full twists along the $c_i$'s we can find a hyperbolic structure $X'$ with the same length parameters, while the twists parameters are between $0$ and $\ell_i(X')=\ell_i(X)$, and $M_P(X)\overset{isom}\cong M_P(X')$.

By doing twists of length at most $\ell_i(X)/4$ we get a surface $X_s$ with the same length parameters and all seams of pair of pants matching up. Then, the twist parameters are equal to either zero or $\ell_i(X)/2$.
\epf

\section{The renormalized volume of symmetric surfaces}
\label{sc:symmetric}

In this section we estimate the renormalized volume of a Schottky filling of a surface $X\in \T_s$ corresponding to a ``symmetric'' pants decomposition. This will be used in the proof of Theorem \ref{mainthm}. In the next two sections, we will bound the variation of the renormalized volume under a variation of the twist parameters in the Fenchel-Nielsen coordinates, and as a consequence we will be able to obtain an upper bound on the renormalized volume of Schottky fillings which are non-symmetric by comparing their renormalized volume to that of a symmetric surface obtained by changing the twist parameters.

 In the following lemma we will deal with manifold whose convex core is 2-dimensional. Thus, it will be useful to use a slightly modified definition of convex core boundary which is more compatible with the corresponding conformal boundary. In what follows we denote by $\partial \widehat C(M)$ the ``boundary'' of $C(M)$ for $M$ any convex co-compact hyperbolic manifold and we define:
  \begin{itemize}
  \item $\partial \widehat C(M)$ is the boundary of $C(M)$ in the usual sense if $C(M)$ has non-empty interior,
  \item if $C(M)$ is a totally geodesic orientable surface $\Sigma \subset M$, then $\partial \widehat C(M)$ is the union of two copies of $\Sigma$ with opposite orientation, if $\partial\Sigma\neq\emp$ then the two copies of $\Sigma$ are glued along their common totally geodesic boundary.
   \item if $C(M)$ is a totally geodesic non-orientable surface $\Sigma \subset M$, then $\partial\widehat C(M)$ is the orientation double-cover of $\Sigma$.
  \end{itemize}
     In all cases, $\partial\widehat C(M)$ is homeomorphic to $\partial_\infty M$. Specifically, the hyperbolic Gauss map, which sends a unit vector normal to a support plane of $C(M)$ to the endpoint at infinity of the geodesic ray it defines, is a homeomorphism from the unit normal bundle of $C(M)$ -- which is itself homeomorphic to $\partial\widehat C(M)$ -- to $\partial_\infty M$.

  The ``boundary'' $\partial\widehat C(M)$ is equipped with an induced metric $m$, which is hyperbolic. However, it is {\em pleated} along a {\em measured lamination} $\beta$ which is geodesic for $m$, with the transverse measure recording the amount of pleating along the leaves, see \cite{thurston-notes,epstein-marden}. When $C(M)$ is a totally geodesic surface $\Sigma$, the support of $\beta$ corresponds to the boundary of $\Sigma$, with each leave equipped with a weight $\pi$.

  Let $X$ be the conformal structure at infinity of $M$. Then $X$ is obtained from $m$ and $\beta$ by a geometric construction called {\em grafting}, see e.g. \cite{dumas-survey}. Given a closed surface $S$ of genus at least $2$, this grafting operation defines a map
\[ \gr:\T (S)\times \cML (S)\to \T (S)~, \]
  where $\cML (S)$ denotes the space of measured laminations on $S$. The key property that is important to us here is a result of Scannell and Wolf \cite{scannell-wolf}: if $\lambda\in\cML (S)$ is fixed, the map $\gr(\cdot,\lambda):\T (S)\to \T (S)$ is a homeomorphism.

Now, if $M$ is a convex co-compact hyperbolic manifold with convex core $C(M)$ a totally geodesic surface $\Sigma$ with boundary, 
then the lamination $\beta$ is the fixed-point set of an orientation-reversing involution $\sigma:\partial\widehat C(M)\to \partial\widehat C(M)$ such that $\Sigma=\partial\widehat C(M)/\sigma$, and the induced metric $m$ on $\partial\widehat C(M)$ is invariant under $\sigma$. Since $m$ and $\beta$ are both invariant under $\sigma$, so is the conformal structure at infinity $X=\gr(m,\beta)$.

\blem \label{flatCC}{
Let $\sigma: S\to S$ be an orientation-reversing involution with $\text{Fix}(\sigma)\neq\emp$ and quotient surface $\Sigma=S/\sigma$. Then, for any invariant conformal structure $X\in \T_\sigma(S)$ there exists a handlebody $H$ with a convex co-compact hyperbolic structure such that the convex core of $H$ is homeomorphic to $\Sigma$ and the conformal boundary of $H$ is $X$.}
\elem

\begin{proof}

 Let $\beta=\partial \Sigma$. We claim that the restriction map
\[ \gr(\cdot, \beta)_{|\T_\sigma(S)}:\T_\sigma(S)\to \T_\sigma(S) \]
  is onto. Indeed, let $X\in \T_\sigma(S)$. Since $\gr(\cdot, \beta):\T (S)\to \T (S)$ is a homeomorphism, there exists a unique $Y\in \T (S)$ such that $\gr(Y,\beta)=X$. But then
\[ \gr(\sigma^*Y,\beta)=\gr(\sigma^*Y,\sigma^*\beta)=\sigma^*\gr(Y,\beta)=\sigma^*X=X=\gr(Y,\beta)~, \]
  and since $Y$ is unique, $\sigma^*Y=Y$, so that $Y\in \T_\sigma(S)$.

{ Let $Y_\sigma=Y/\sigma$, homeomorphic to $\Sigma$, be the quotient surface of the hyperbolic surface $Y$ by the locally isometric involution $\sigma$. Then, $Y_\sigma$ has a uniformization $\Gamma < \text{Isom}^\pm(\mathbb H^2)$}. By considering $\Gamma$ inside $\text{Isom}^\pm(\mathbb H^3)$, by the natural inclusion, and the corresponding quotient $\mathbb H^3/\Gamma$, we obtain a hyperbolic 3-manifold whose convex core is $Y_\sigma$. By the above discussion, we also know that the conformal boundary is $X$. The fact that $H=\mathbb H^3/\Gamma$ is a handlebody follows from the fact that $H$ is homeomorphic to either $Y_\sigma\times I$, if $Y_\sigma$ has non-empty boundary or is orientable, or to the twisted bundle $Y_\sigma\overset\sim\times I$, if $Y_\sigma$ is non-orientable with empty boundary. As $\Sigma\cong Y_\sigma$ has boundary by Remark \ref{generalhandlebody} this yields a handlebody.  \end{proof}

\brem { The following remark is not needed in the rest of this paper, however, one should note that Lemma \ref{flatCC} also works for fixed point free involution if one allows for the topological condition to be that of a twisted $I$-bundle $K\overset\sim\times I$ over a closed non-orientable surface $K$.}
\erem
\brem\label{remallloopsfixed} In the case we have a pants decomposition $P$ such that, for each $\gamma\in P$, $\sigma(\gamma)=\gamma$, we can also infer from Lemma \ref{flatCC} and Lemma \ref{symmquotient} that $H\in CC_P(S)$, i.e. $P$ compresses in $H$ and the seams of $P$ form $\partial X_\sigma$.
\erem 

For a convex co-compact hyperbolic 3-manifold $M$, 
\begin{equation}\label{boundVR} V_R(M)\leq V_C(M)-\frac 1 4 L(\beta_M)~,\end{equation}
see \cite[Lemma 4.1]{compare} (and also \cite[Theorem 3.7]{BBB2018}). In the case considered here, the bending lamination is given by a multi-curve with bending measure given by assigning the weight $\pi$ to each curve, see Lemma \ref{flatCC}. Then, its length is given by:
\begin{equation}\label{boundVR2} L(\beta_M)=\pi\sum_{\gamma\in\pi_0}\ell_Y(\gamma)~,\end{equation}
for $Y$ the hyperbolic structure on the convex-core boundary and $\pi_0$ the set of the simple closed curves composing the multicurve. Thus, one has $L(\beta_M)> 0$ and so, by Lemmas \ref{symmquotient} and \ref{flatCC} we obtain the following statement.

\bthm
Let $X\in \T_s(S)$, and let $\sigma:S\rar S$ be such that $X\in \mathcal T_\sigma(S)$. Then there exists a handlebody filling $H_X$ such that 
\[V_R(H_X)\leq -\f \pi 4 \ell_{X_\sigma}(\partial X_\sigma)<0\ .\]
\ethm 
\bpf {  By equation \eqref{boundVR} and \eqref{boundVR2}, as the convex core volume of $H_X$ is zero, we have:
\[V_R(H_X)\leq -\f \pi 4 \ell_{Y_\sigma}(\partial \Sigma)~,\]
for $Y_\sigma$ the hyperbolic structure on the convex core induced by Lemma \ref{flatCC}. Recall that by Lemma \ref{symmquotient} and Remark \ref{remallloopsfixed}, as isotopy classes of loops in $S$, we have that $\partial \Sigma$ and $\partial X_\sigma$ are the same. We remark that, with respect to the metric obtained by grafting $Y_{\sigma}$ along $\partial\Sigma$ (without reuniformizing), i.e., along the bending lamination of $H_X$, the boundary $\partial\Sigma$ has the same length as in $Y_{\sigma}$ (see \cite[Section 4.1]{dumas-survey} for the geometric definition of grafting). Moreover, the \textit{grafting metric} coincides with the \textit{Thurston metric}, which is conformal to $X$, and defined as \[h_{Th}(z)=\inf_{\Omega(\Gamma)}h_D(z)~,\]
where $H_X=\mathbb{H}^3/\Gamma$, and the infimum is taken on the round disks $D$ immersed in $\Omega(\Gamma)$, with $h_{D}$ the hyperbolic metric on $D$.
By the Schwarz Lemma, the Thurston metric is bigger then the hyperbolic metric at infinity. Therefore
\[  V_R(H_X)\leq -\f \pi 4 \ell_{Y_\sigma}(\partial \Sigma)\leq -\f \pi 4 \ell_{X_\sigma}(\partial X_\sigma).\] } 
\epf
If we know some curves are short in $X$ and the pants decomposition is fixed component by component we obtain the following estimate, Lemma \ref{Wvolsym}. This is the main such estimate we will use in this work. For completeness we also prove the other option in Lemma \ref{Wvolsym2}.

\blem \label{Wvolsym}
 There exist universal constants $S, Q\geq 0$ as follows. Let $X\in\mathcal T_s(S)$, $\sigma$ so that $X\in \T_\sigma(S)$, $M=M_P(X)$ be the Schottky manifold corresponding to any pants decomposition for which each curve is invariant under {$\sigma:S\rar S$}, and such that there are $0\leq k\leq 3g-3$ geodesic loops of $P$ of length $\ell_X(\gamma_i)\leq 1$. Then,
\[ V_R(M_P(X))\leq -\f S4\sum_{i=1}^k \frac 1 {\ell_X(\gamma_i)} +\f Q 4 k\leq \f k 4 \left( - S +Q\right)<0~. \] 
Specifically, one can take $S=\frac{4\pi^3}{\sqrt e}\sim75.225$ and $Q=4\pi \log\left(\frac{\pi e^{0.502\pi}}{\argsinh(1)} \right)\sim 35.7901\leq 36$.
\elem

\bpf
By Lemma \ref{flatCC} the convex-core of $M$ is a totally geodesic surface {\sout{isometric to $X_\sigma$}} and so $V_C(M)=0$. However, by Theorem 2' of \cite{BC2005}, we have
\[ \sum_{i=1}^k \left( \frac S {\ell_X(\gamma_i)} -Q\right) \leq L(\beta_M). \]
Applying it to equation \eqref{boundVR}, one gets:
\[  V_R(M_P(X))\leq -\frac 14\sum_{i=1}^k\left( \frac S {\ell_X(\gamma_i)} -Q\right)=-\f S4\sum_{i=1}^k \frac 1 {\ell_X(\gamma_i)} +\f Q 4 k\leq \f k 4 \left( -S +Q\right), \] 
concluding the proof.\epf
 
The case in which the pants curves are not fixed component-wise requires introducing some auxiliary functions from \cite[Corollary 1]{BC2003}, these functions will only be needed here. For $m=\cosh^{-1}(e^2)$ we define
$$ g(x)=e^{-m}\frac {e^{-\pi^2/2x} }{2} $$
and
$$ F(x)=\frac x 2+\sinh^{-1}\left(\frac{\sinh(x/2)}{\sqrt{1-\sinh^2(x/2)}}\right)~. $$
Since $F$ is invertible we let $K(x)=\frac{2\pi}{F^{-1}(x)}$ and then define $L(x)=1+K(g(x))$.

\blem \label{Wvolsym2}   Let $X\in\T_s(S)$, and $M(X)$ be the Schottky manifold with {flat convex core and conformal boundary $X$}. Let $\mathfrak m=\set{\gamma_1,\dotsc,\gamma_k}$ be the collection of geodesic loops point-wise invariant by $\sigma$ and let $\rho_X$ be half of the length of the shortest simple closed compressible geodesic in $X$. Then, \[ V_R(M(X))\leq -\f {\pi}{4L(\rho_X)}\sum_{i=1}^k \ell_X(\gamma_i) ~.\]
\elem

\bpf
By Lemma \ref{flatCC} the convex-core of $M$ is a totally geodesic surface {$Y$\sout{isometric to $X_\sigma$}} and so $V_C(M)=0$. Moreover, by Lemma \ref{symmquotient} $\partial X_\sigma$ is given by the multi-curve $\mathfrak m$ of geodesic loops that are point-wise fixed by $\sigma$. Thus, by Corollary 1 of \cite{BC2003} we have:
\[ \ell_{Y}(\gamma)\geq \f 1 {L(\rho_X)} \ell_X(\gamma)~.\]
Then, by applying it to equation \eqref{boundVR} we obtain the required result. \epf

\section{Variation of the renormalized volume under an earthquake}
\label{sc:variation}

In this section we compute how the renormalized volume changes under earthquake paths in the deformation space. 
\subsection{First-order variation of the renormalized volume}
\label{ssc:dehn-inf}

We start the section with a formula for $\diff V_R$ at $M_P(X)=\mathbb{H}^3/\Gamma$. Recall that by $\mathcal S(f)$ we are denoting the Schwarzian derivative of the developing map of the domain of discontinuity $\Omega(\Gamma)$ of the Schottky hyperbolic $3$-manifold $M_P(X)$, and by $S$ the boundary $\partial\bar M$. {We will sometimes refer to $\mathcal S(f)$ just as the Schwarzian of $M_P(X)$}.

\blem\label{intbypart}
  Let $\mu$ be an infinitesimal earthquake (at unit velocity) along a simple closed geodesic on $X$, parameterised at unit velocity by $\gamma:\R/ \ell\Z \to X$. Then, for $q=\mathcal S(f)$:
  
\[ \diff V_R(\mu)=-\f 12\int_{\R/\ell\Z}\re(q(i\dot\gamma(t),\dot\gamma(t))) dt=\im\left(\f 12\int_{\R/\ell\Z} q(\dot\gamma(t),\dot\gamma(t))\right)~.\]
\elem

\bpf
  Let $v$ be a vector field realizing the infinitesimal earthquake along the image of $\gamma$. That is, $v$ is the vector field obtained by differentiating at zero, with respect to the time parameter $t$, the family of diffeomorphisms $\phi_{\gamma}(t)$ corresponding to a length $t$ earthquake along $\gamma$. We assume that $v$ is discontinuous along $\gamma(\R/\ell\Z)$, that is, it has limit zero on the right side and equal to $\dot\gamma(t)$ along $\gamma(\R/\ell\Z)$ { and is continuous on the left}.

  The first-order variation of the complex structure associated to $v$ is then determined by the Beltrami differential $\mu=\bar\partial v$. {Where here by $\bar\partial v$ we mean the $L^\infty$ weak$^*$ limit of $\bar\partial v_n$ for $v_n$ smooth compactly supported vector fields that are $C^\infty$ approximations of $v$, converging in the uniform topology on compact sets of $S\setminus \gamma$.} Specifically, choosing a complex coordinate $z$, we can write:
\[ v=2\omega \left(\partial_z+\partial_{\bar z}\right)~, \]
and note that $\omega$ vanishes on the right half-neighbourhood of $\gamma$.

Consider the area form $dx\wedge dy$ associated to $z=x+iy$, and note that $d\bar{z}\wedge dz = 2i (dx\wedge dy)$. We have
\[\bar\partial v = 2(\bar\partial\omega) (d\bar z\otimes \partial_z+d\bar z\otimes \partial_{\bar z})~, \]
and so if $q=g(z)dz^2$, 
\begin{align*}
 \langle q,\bar\partial v\rangle &=  \f{1}{2i} \int_{S}2g(z)(\bar\partial \omega(z) )(d\bar z\wedge( dz^2(\partial_z))+d\bar z\wedge (dz^2(\partial_{\bar z})))\\
 &=\int_{S} 2g(z)(\bar\partial w(z))dx\wedge dy
 \end{align*}
  by definition of the duality product commonly used between Beltrami differentials and holomorphic quadratic differentials, see \cite{Hubbard2016}.

  Consider now the one-form defined by $\alpha=q(v,\cdot)=2\omega(z) g(z) dz$. Then
  \begin{align*}
    \bar\partial \alpha 
               & = \bar\partial (2\omega(z) g(z) dz) \\
               & =  2(\bar\partial\omega(z))g(z) d\bar z \wedge d z+2\omega(z) \bar\partial g(z) d\bar z\wedge dz\\
               & =  4i (\bar\partial \omega(z)) g(z)dx\wedge dy~,
  \end{align*}
  because $g$ is holomorphic, $\bar\partial g=0$, and $d\bar{z}\wedge dz =2i (dx\wedge dy) $.

  The outcome of this discussion is that
  \[ \int_S\bar\partial \alpha =\int_S 4i g(z) (\bar\partial \omega(z) ) dx\wedge dy=2i\int_S 2g(z)(\bar\partial \omega(z)) dx\wedge dy=2i \langle q,\bar\partial v\rangle~.\]
  
  Therefore, we get
\[ \langle q,\mu\rangle =- \frac i2\int_S \bar\partial\alpha~. \]
  However, $\alpha$ is a complex 1-form, so that $\partial \alpha=\partial(2gw) dz\wedge dz=0$, and as a consequence 
\[ d\alpha = (\partial + \bar\partial) \alpha = \bar\partial \alpha~. \]
  Using Stokes on $S'=S\setminus \gamma(\R/\ell\Z)$, we obtain that, since $\alpha$ vanishes on one component of $\partial S'$:
\[\langle q,\bar \partial v\rangle =- \frac i2\int_{S} d\alpha =- \frac i2\int_{\partial S'} \alpha(\dot\gamma(t)) dt= -\f i2\int_0^\ell \alpha(\dot \gamma (t))dt~. \]
 However, by definition of $\alpha$ we obtain that
 \[  \langle q,\bar \partial v\rangle =-\f i2\int_0^\ell q(v\vert_{\gamma(t)},\dot\gamma(t))dt=-\f 12\int_0^\ell q(i\dot\gamma(t),\dot\gamma(t))dt~.\]
  
  The first order variation of the renormalized volume, thanks to Theorem \ref{dVRdiff}, is equal to:
  \begin{align*} \diff V_R(\mu)&= \re\left(\langle q,\bar \partial v\rangle \right)=\re\left( -\f 12\int_0^\ell q(i\dot\gamma(t),\dot\gamma(t))dt\right)\\
                          &=-\f 12\int_0^\ell \re\left( q(i\dot\gamma(t),\dot\gamma(t))\right)dt~,
  \end{align*}
completing the proof.
\epf

\bdefi An {\em earthquake path} $c_{\mathbf t} :[0,1]\rar CC_P(S)$, with $\mathbf t=(t_1,\dotsc, t_{3g-3})$, {is a path which at time $s\in[0,1]$ twists $st_i\in \R$} along each pants curve $\gamma_i\in P$ {of $c_{\mathbf t}(0)$}. 
\edefi

 For a loop $\gamma$ we use $\operatorname{inj}\vert_{\gamma}$ to denote {half of} the length of the shortest { loop $\delta$ such that  $\delta$ is compressible and $ \delta$ is either $\gamma$ or $\delta$ intersects $\gamma$ essentially in $M$}. Note that if $\gamma$ is a compressible geodesic loop of length $\leq\epsilon_0$ in $X$, then $\operatorname{inj}\vert_{\gamma}=\f{\ell_X(\gamma)}2$.

\blem \label{dVrexp}{ Let $c_t(s)$, for $s\in[0,1]$ and a fixed $t\in \mathbb{R}$, be an earthquake path along a simple geodesic loop $\gamma$ starting at the Riemann surface $X_0$. Then, the following bound for $\abs{\diff (V_R\circ c_t)}$ holds at any $s\in[0,1]$:
\[ |\diff (V_R\circ c_t)|\leq 3\ell_{X_0}(\gamma)\coth^2\left(\frac {\operatorname{inj}\vert_{\gamma}} 2\right)t~.\]
In particular, if $\operatorname{inj}\vert_{\gamma}\geq 1/2$ we have 
\[ |\diff (V_R\circ c_t)|\leq 3\ell_{X_0}(\gamma)\coth^2\left(\frac {1} 4\right)t=C\ell_{X_0}(\gamma)t~, \quad C=3\coth^2\left(\frac {1} 4\right)<50.013~.\]}
\elem

\begin{proof} {First, observe that the length of $\gamma$  remains constantly equal to $\ell_{X_0}(\gamma)$ along the earthquake path $c_t(s)$. Moreover, since earthquaking forms a flow (i.e. $c_t(s_1+s_2)=c_t(s_1)\circ c_t(s_2)$), the scaling by $t$ of the infinitesimal earthquake $\mu_s$ along $\gamma$ at $X_s=c_t(s)$ coincides with the derivative of $c_t(s)$ at $s$. Then, at any $s\in[0,1]$, we can use the integration by part of Lemma \ref{intbypart}. Denoting by $\mathcal S(f_s)(z)=q_s(z)dz^2$ the Schwarzian associated through uniformization to $c_t(s)$, we can estimate $\abs{ q_s(z)} \leq 6\coth ^2\left(\frac {\operatorname{inj}\vert_{\gamma}} 2\right)$ (see \cite[Corollary 2.12]{BBB2018}, and note that the factor $4$ comes from the hyperbolic metric), yielding the first bound. The second estimate follows by direct computation.} \end{proof}

\subsection{Earthquake paths and $V_R$ estimates}
\label{ssc:dehn}

In this section we compute the change of renormalized volume under a path $c_{\mathbf t}:[0,1]\rar CC_P(S)$ obtained by doing earthquakes along geodesic loops in the pants decomposition $P$.

\bthm\label{dehntwistpath}
Let $c_{\mathbf t} :[0,1]\rar CC_P(S)$ be an earthquake path, and let $\ell_i=\ell_{X_0}(\gamma_i)$. Then 
\[ \abs{ V_R(X_1) -V_R(X_0) }\leq   \sum_{i=1}^{k} (3\ell_i\coth^2 \left( \ell_i/4\right) )t_i +C \sum_{i=k+1}^{3g-3} t_i\ell_i~,\]
where $\gamma_1,\dotsc, \gamma_k$ are the geodesic loops of $P$ with $\ell_i<1$ and for all $j>k$ we have $2\operatorname{inj}\vert_{\gamma_j}\geq 1$, {and $C=3\coth^2(1/4)$}.
\ethm

\bpf Pick a $1$-thick/thin pants decomposition with $k$ geodesic loops less than $1$ and integrate Lemma \ref{dVrexp}.
\epf

Since, by Lemma \ref{FNtwist}, to reach a symmetric surface we need to twist at most $\ell_{X}(\gamma_i)/4$, we can take $t_i\leq \ell_{X}(\gamma_i)/4$ in the above expression and obtain the following statement.

\bcor\label{mainest}
Let $X\in\T(S)$ and $P=\set{\gamma_i}_{i=1}^{3g-3}$ be a pants decomposition in which the first $k$ curves have length less than $1$ and the others have injectivity radius at least $1$. Then, there exists a symmetric surface $X_0$ such that 

\[\abs{ V_R(X) -V_R(X_0) }\leq  \f 34\sum_{i=1}^{k} \coth^2 \left( \ell_i/4\right) \ell_i^2+\f C4\sum_{i=k+1}^{3g-3}\ell_i^2,\]
with $\ell_i=\ell_{X_0}(\gamma_i)$ and $C=3\coth^2\left(\frac {1} 4\right)<50.013$.
\ecor

The above estimates also work in the setting of general convex co-compact manifolds. Let $CC(M)$ be the deformation space which is also parameterised by the quotient of $\mathcal T(\partial M)$ by Dehn twists along disks. Let $c_{\mathbf t}^{\mathfrak m}:[0,1]\rar CC(M)$ be an earthquake path along a multi-curve $\mathfrak m\subset S$. Assume that the multi-curve $\mathfrak m$ can be subdivided, according to the reference metric $X_0$, in the following way:
\begin{itemize}
\item $\mathfrak m_{1}^c$ is the set of geodesic loops $\gamma$ of $\mathfrak m$ that are compressible and have length at most $1$;
\item $\mathfrak m_{ 1}$ is the set of geodesic loops $\gamma$ contained in compressible components of $\partial M$ and not in $\mathfrak m_1^c$, and such that any compressible loop intersecting $\gamma$ essentially has length at least $1$;
\item $\mathfrak m_\infty$ is the set of geodesic loops $\gamma$ of $\mathfrak m$ that are contained in components of $\partial M$ that are incompressible.
\end{itemize}
Note that not every $\m$ admits such a decomposition with respect to the given $X_0$.

\bthm\label{dehntwistpath2} Let $X_0\in\T(\partial M)$ and $\m=\m_{1} ^c\cup \mathfrak m_{1}\cup \m_\infty$ be a multi-curve and $c_{\mathbf t}^{\mathfrak m} $ be an earthquake path terminating at $X_1$. Then 
\[ \abs{ V_R(X_1) -V_R(X_0) }\leq   \sum_{\gamma_i\in \pi_0(\mathfrak m_{1}^c)} (3\ell_i\coth^2 \left( \ell_i/4\right) ) t_i +C \sum_{\alpha_j\in\pi_0(\mathfrak m_{ 1} )} t_j\ell_j +{3}\sum_{\beta_k\in\pi_0(\mathfrak m_\infty)}  t_k \ell_k~,\] 
for  $C=3\coth^2\left(\frac {1} 4\right)<50.013$.

\ethm
\bpf The first two cases follow by the previous computations and integrating Lemma \ref{dVrexp}. For the last case we bound the norm of the Schwarzian on the geodesic loops in $\mathfrak m_\infty$ by the Kraus-Nehari estimate \cite{Ne1949,kraus1932zusammenhang} and then integrating gives the result.
\epf

  \section{Main Results}
  \label{sc:proofs}

We now put together the results from the previous sections to prove the main Theorem \ref{mainthm} and Corollary \ref{maincoro}. 

\begin{customthm}{\ref{mainthm}}
Let $X$ be a closed Riemann surface of genus $g\geq 2$. Assume that there are $k$ disjoint simple closed curves $\gamma_1, \cdots, \gamma_k$ such that $\ell(\gamma_i)\leq 1, 1\leq i\leq k,$ and there are no other geodesic loops of length less or equal to $1$ in $X$. Then there exists a pants decomposition $P$ containing the $\gamma_i$'s such that
\[ V_R(M_P(X))\leq -\f{\pi^3}{\sqrt e}\sum_{i=1}^k \frac 1{\ell(\gamma_i)} + \left(9+\frac 3 4\coth^2\left(\f14\right)\right)k + 81\coth^2\left(\f14\right)\pi(3g-3-k)(g-1)^{2}~. \] 
\end{customthm}
 \bpf
 Let $ P$ be a pants decomposition containing the $k$ geodesic loops $\gamma_1,\dotsc,\gamma_{k}$ shorter than $1$ and the $\alpha_i$, $i=k+1,\dotsc, 3g-3$, being {\em Bers pants curves} (see \cite[Theorem 12.8]{FM2011}). 

That is, we have:
\begin{itemize}
\item $\ell_X(\gamma_i)\leq 1$ for $i\leq k$;
\item $1<\ell_X(\alpha_i)\leq B_g\leq 6\sqrt{3\pi}(g-1)$, see \cite[Theorem 5.1.4]{Bu1992}, and $\operatorname{inj}\vert_{\alpha_i}\geq 1$ for $k<i\leq 3g-3$;
\item $ P$ has seams such that in the $FN$ coordinates induced by $ P$, $FN(X)$ has no twists bigger than $\ell_X(\gamma_i)/4$ or $\ell_X(\alpha_i)/4$ (see Lemma \ref{FN} and Lemma \ref{FNtwist}).
\end{itemize}
Let $c_{\mathbf t}$ be the path in $FN$ coordinates from $X$ to $X_s$, the symmetric surface. Then, $c_{\mathbf t}$ can be thought of doing $3g-3$ twists along each pants curve, each of length at most $\ell_X(\gamma_i)/4$ or $\ell_X(\alpha_i)/4$, see Lemma \ref{FNtwist}. Then, for $C=3\coth^2\left(\frac {1} 4\right)$, by Corollary \ref{mainest} we get:
\begin{align*}
 \abs{V_R(X)-V_R(X_s)} &\leq  \f 34\sum_{i=1}^{k} \coth^2 \left( \ell_i/4\right) \ell_i^2+\f C4\sum_{i=k+1}^{3g-3}\ell_i^2  \\
&\leq  \frac C 4k+\frac C 4 \sum_{i=k+1}^{3g-3}B_g^2\\
&\leq \frac C 4 k +\frac C 4  (3g-3-k)B_g^2 \\
&\leq  \frac C 4k+27C \pi(3g-3-k)(g-1)^2, \end{align*} 
where we used the fact that $B_g\leq  6\sqrt{3\pi}(g-1)$ and $\coth^2 \left( x/4\right) x^2$ is an increasing function. Thus, we get that:
\[ V_R(X)\leq V_R(X_s) +\frac C 4k+27 C\pi(3g-3-k)(g-1)^2.\]
Since $\ell_i\leq1$ for $i\leq k$ by using Lemma \ref{Wvolsym} to estimate $V_R(X_s)$ we have:
 \[V_R(X_s)\leq -\frac 14\sum_{i=1}^{k}  \left( \frac S {\ell_X(\gamma_i)} -Q\right), \]
 for $S=\frac{4\pi^3}{\sqrt e}$ and $Q=4\pi \log\left(\frac{\pi e^{0.502\pi}}{\argsinh(1)} \right)\sim 35.7901\leq 36$. Then, we obtain the following bound:
 
 \begin{align*} 
 V_R(X)&\leq \sum_{i=1}^{k}\left(- \frac S {4\ell_X(\gamma_i)}+\f Q4\right)+\frac C 4k+27C\pi(3g-3-k)(g-1)^2\\
 &\leq \sum_{i=1}^{k}\left(- \frac {\pi^3} {\sqrt e\ell_X(\gamma_i)}\right)+9k+\frac C 4k+27C\pi(3g-3-k)(g-1)^2\\
 &\leq -\frac{\pi^3}{\sqrt e} \sum_{i=1}^{k}\left(\frac{1}{\ell_X(\gamma_i)}\right)  +\left (9+\frac C 4\right)k+27C\pi(3g-3-k)(g-1)^2~.
 \end{align*}
 Substituting for $C=3\coth^2\left(\frac {1} 4\right)$ concludes the proof.
 \epf

\begin{customcor}{\ref{maincoro}} For all $g\in \N$ s.t. $g\geq 2$, $0<k\leq 3g-3$ and $0<k_1\leq k$ there exists a constant $A=A(g,k_1,k-k_1)>0$ such, that if $X$ is a Riemann surface with $k_1$ geodesic loops of length less than $A$ and $k$ geodesic loops of length at most $1$, then $X$ admits a Schottky filling with negative renormalized volume. \end{customcor}

 \bpf Let $ P$ be a pants decomposition containing the $k_1$ geodesic loops, $\gamma_1,\dotsc,\gamma_{k}$ shorter than $A$ and $k-k_1$ geodesic loops $\gamma_{k_1+1},\dotsc,\gamma_k$ of length at most $1$ and the $\alpha_i$, $i=k+1,\dotsc, 3g-3$ are Bers pants curves. 

That is, we have:
\begin{itemize}
\item $\ell_X(\gamma_i)< A$ for $1\leq i\leq k_1$;
\item $\ell_X(\gamma_i)\leq 1$ for $k_1<i\leq k$;
\item $1<\ell_X(\alpha_i)\leq B_g\leq 6\sqrt{3\pi}(g-1)$ and $\operatorname{ inj}\vert_{\alpha_i}\geq 1$ for $k<i\leq 3g-3$.
\end{itemize}
Then, by Theorem \ref{mainthm} we get:
\[ V_R(M_P(X))\leq -\f{\pi^3}{\sqrt e}\sum_{i=1}^k \frac 1{\ell(\gamma_i)} +\left (9+\frac C 4\right)k  + 27C\pi(3g-3-k)(g-1)^{2}~, \]
which can be further decomposed in:
\[ V_R(M_P(X))\leq  -\f{\pi^3}{\sqrt e}\left(\sum_{i=1}^{k_1} \frac 1{\ell(\gamma_i)}+\sum_{i=k_1+1}^{k} \frac 1{\ell(\gamma_i)}\right)  +\left (9+\frac C 4\right)k + 27C\pi(3g-3-k)(g-1)^{2}~. \]
Since for $i\leq k_1$ we have that $\frac {1}{\ell_X(\gamma_i)}\geq \frac 1 A$ and, similarly, for $k_1+1\leq i\leq k$ we have that $\frac {1}{\ell_X(\gamma_i)}\geq 1$ we get:
 \[ V_R(X)\leq -\f{\pi^3}{\sqrt e}\left(\frac {k_1} A + k-k_1\right)+ \left (9+\frac C 4\right)k+ 27C\pi(3g-3-k)(g-1)^{2}~. \] 
We want to find an upper bound on $A$ that makes the above expression negative. Note that
\[ B\eqdef -\f{\pi^3}{\sqrt e}(k-k_1)+\left (9+\frac C 4\right)k + 27C\pi(3g-3-k)(g-1)^{2}> 2k>0 ~,\] 
as the smallest case for $B$ is for $k=3g-3$ and $k_1=0$. Then, to have
\[  -\f{\pi^3}{\sqrt e}\frac {k_1} A +B<0~,\]
 it suffices to take:
 \[ A< \f{\pi^3}{\sqrt e}\f{k_1} {B},\]
 concluding the proof.\epf

\textbf{Conflict of interests}: No conflicts of interests to declare.
\thispagestyle{empty}
{\small
\markboth{References}{References}
\bibliographystyle{abbrv}
\bibliography{mybib_3}{}
}

\end{document}

%% file: packages.tex
\usepackage{latexsym,enumitem}
\usepackage{amssymb}
\usepackage[cp850]{inputenc}
\usepackage{epsfig}
\usepackage{psfrag}
\usepackage{amsthm}
\usepackage{amscd}
\usepackage{amsmath}
\usepackage{amsfonts}
\usepackage{graphics,caption}
\usepackage[all]{xy}
\usepackage{etoolbox}
\patchcmd{\quote}{\rightmargin}{\leftmargin 2em \rightmargin}{}{}
 \usepackage[sort,numbers]{natbib} 
 \usepackage{mathtools}
\usepackage[normalem]{ulem}
\usepackage{array}

\usepackage[english]{babel}
\usepackage{comment}
\usepackage{color}
\usepackage[colorlinks]{hyperref}

\let\phi\varphi

\DeclareMathOperator{\eqdef}{\coloneqq} %definition
\newcommand{\f}[2]{\frac{#1}{#2}} %Fractions

\let\epsilon\varepsilon
\let\subset\subseteq
\newcommand{\be}{\begin{equation*}}
 \newcommand{\ee}{\end{equation*}}
\newcommand{\bpf}{\begin{dimo}}
\newcommand{\epf}{\end{dimo}}
\newcommand{\bdefi}{\begin{defin}}
\newcommand{\edefi}{\end{defin}}
\newcommand{\bthm}{\begin{thm}}
\newcommand{\ethm}{\end{thm}}
\newcommand{\blem}{\begin{lem}}
\newcommand{\elem}{\end{lem}}
\newcommand{\bcor}{\begin{cor}}
\newcommand{\ecor}{\end{cor}}

\newcommand{\bprop}{\begin{prop}}
\newcommand{\eprop}{\end{prop}}
\newcommand{\bese}{\begin{ese}}
\newcommand{\eese}{\end{ese}}
\newcommand{\brem}{\begin{rem}}
\newcommand{\erem}{\end{rem}}
\newcommand{\bpfc}{\begin{dimoclaim}}
\newcommand{\epfc}{\end{dimoclaim}}

\newcommand{\rar}{\rightarrow} %arrow right
\newcommand{\diff}{\mathop{}\!\mathrm{d}} %exact differential
	
\newcommand{\abs}[1]{\left\lvert#1\right\rvert}						%mod
\newcommand{\set}[1]{\left\{#1\right\}}					%set, curly brackets
\newcommand{\T}[1]{\mathsf{{T_{#1}}}}				%Separation Axioms
\newcommand{\quotient}[2]{\left.\raisebox{.1em}{$#1\!$}\middle/\raisebox{-.1em}{$#2$}\right.}

\DeclareMathOperator{\emp}{\varnothing} %Known Sets, Fields and so on...
\DeclareMathOperator{\N}{\mathbb N}			%Naturals (with 0)
\DeclareMathOperator{\R}{\mathbb R}			%Reals
\DeclareMathOperator{\C}{\mathbb C}			%Complexes
\DeclareMathOperator{\Z}{\mathbb Z}			%Integers
\newcommand{\m}{\mathfrak m}

\DeclareMathOperator{\im}{Im} %image
\DeclareMathOperator{\id}{id} %identity map

\newcommand{\M}{\mathscr M}

\newenvironment{quot}
{
	\vspace{-0.2cm}
	\vspace{0.2cm}
}

\theoremstyle{definition}
\newtheorem{d1}{Definition}[section] %fa partire tutto da [section]

\newenvironment{defin}
{
	\begin{quot}
		\begin{d1}
		}
		{\end{d1}
	\end{quot}

}

\theoremstyle{definition}
\newtheorem{r1}[d1]{Remark}%[section]

\newenvironment{rem}
{
	\begin{quot}
		\begin{r1}
		}
		{\end{r1}
	\end{quot}
}

\theoremstyle{definition}
\newtheorem{e1}[d1]{Exercise}%[section]

\theoremstyle{definition}
\newtheorem{ese1}[d1]{Example}

\newenvironment{ese}
{
	\begin{quot}
		\begin{ese1}
	}
	{	
		\end{ese1}
	\end{quot}
}

\theoremstyle{definition}
\newtheorem{f1}[d1]{Formulation}

\theoremstyle{definition}
\newtheorem{f2}[d1]{Fact}

\theoremstyle{definition}
\newtheorem{question}[d1]{Question}

\theoremstyle{definition}
\newtheorem{conjecture}[d1]{Conjecture}

\theoremstyle{definition}
\newtheorem{t1}[d1]{Theorem}%[section]

\newenvironment{thm}
{
	\begin{quot}
		\begin{t1}}
		{\end{t1}
	\end{quot}
}

\theoremstyle{definition}
\newtheorem*{T1*}{Theorem}%[section]

\newenvironment{teor*}
{
	\begin{quot}
		\begin{T1*}}
		{\end{T1*}
	\end{quot}
}

\newenvironment{dimo}
{\begin{proof}[Proof]
	}
	{\end{proof}}

\newenvironment{dimoclaim}{\emph{Proof of Claim:}\;}{\hfill$\square$}

	\theoremstyle{definition}
	\newtheorem{l1}[d1]{Lemma}%[section]
	
	\newenvironment{lem}
	{
		\begin{quot}
			\begin{l1}}
			{\end{l1}
		\end{quot}
	}
	\theoremstyle{definition}
	\newtheorem{p1}[d1]{Proposition}%[section]
	
	\newenvironment{prop}
	{
		\begin{quot}
			\begin{p1}}
			{\end{p1}
		\end{quot}
	}
	
	\theoremstyle{definition}
	\newtheorem{c1}[d1]{Corollary}%[section]
	
	\newenvironment{cor}
	{
		\begin{quot}
			\begin{c1}}
			{\end{c1}
		\end{quot}
	}

\newtheorem{innercustomthm}{Theorem}
\newenvironment{customthm}[1]
  {\renewcommand\theinnercustomthm{#1}\innercustomthm}
  {\endinnercustomthm}

 \newtheorem*{Theorem*}{Theorem}
 \newtheorem*{Proposition*}{Proposition}
 \newtheorem*{Lemma*}{Lemma}

\newtheorem{innercustomcor}{Corollary}
\newenvironment{customcor}[1]
  {\renewcommand\theinnercustomcor{#1}\innercustomcor}
  {\endinnercustomcor}

%% file: thintubes.bbl
\begin{thebibliography}{10}

\bibitem{BBB2018}
M.~Bridgeman, J.~Brock, and K.~Bromberg.
\newblock Schwarzian derivatives, projective structures, and the
  weil--petersson gradient flow for renormalized volume.
\newblock {\em Duke Mathematical Journal}, 168(5):867 -- 896, 2019.

\bibitem{bridgeman-brock-bromberg}
M.~Bridgeman, J.~Brock, and K.~Bromberg.
\newblock Schwarzian derivatives, projective structures, and the
  {W}eil-{P}etersson gradient flow for renormalized volume.
\newblock {\em Duke Math. J.}, 168(5):867--896, 2019.

\bibitem{bridgeman-brock-bromberg:gradient}
M.~Bridgeman, J.~Brock, and K.~Bromberg.
\newblock The {W}eil-{P}etersson gradient flow of renormalized volume and
  3-dimensional convex cores.
\newblock {\em Geom. Topol.}, 27(8):3183--3228, 2023.

\bibitem{bridgeman-bromberg-pallete:convergence}
M.~Bridgeman, K.~Bromberg, and F.~Vargas~Pallete.
\newblock Convergence of the gradient flow of renormalized volume to convex
  cores with totally geodesic boundary.
\newblock {\em Compos. Math.}, 159(4):830--859, 2023.

\bibitem{BC2005}
M.~Bridgeman and R.~D. Canary.
\newblock Bounding the bending of a hyperbolic 3-manifold.
\newblock {\em {Pacific Journal of Mathematics}}, 218(2):299--314, 2005.

\bibitem{BC2003}
M.~J. Bridgeman and R.~Canary.
\newblock From the boundary of the convex core to the conformal boundary.
\newblock {\em Geometriae Dedicata}, 96:211--240, 2003.

\bibitem{brock-bromberg:inflexibility2}
J.~F. Brock and K.~W. Bromberg.
\newblock Inflexibility, {W}eil-{P}eterson distance, and volumes of fibered
  3-manifolds.
\newblock {\em Math. Res. Lett.}, 23(3):649--674, 2016.

\bibitem{Bu1992}
P.~Buser.
\newblock {\em Geometry and spectra of compact {R}iemann surfaces}.
\newblock Modern Birkh\"{a}user Classics. Birkh\"{a}user Boston, Ltd., Boston,
  MA, 2010.
\newblock Reprint of the 1992 edition.

\bibitem{CEM2006}
R.~D. Canary, D.~Epstein, and A.~Marden.
\newblock {\em {Fundamental of Hyperbolic Manifolds}}.
\newblock {Cambridge University Press}, 2006.

\bibitem{skenderis-solod}
S.~de~Haro, K.~Skenderis, and S.~N. Solodukhin.
\newblock Holographic reconstruction of spacetime and renormalization in the
  {A}d{S}/{CFT} correspondence.
\newblock {\em Comm. Math. Phys.}, 217(3):595--622, 2001.

\bibitem{dumas-survey}
D.~Dumas.
\newblock Complex projective structures.
\newblock In {\em Handbook of {T}eichm\"uller theory. {V}ol. {II}}, volume~13
  of {\em IRMA Lect. Math. Theor. Phys.}, pages 455--508. Eur. Math. Soc.,
  Z\"urich, 2008.

\bibitem{EpsteinSurf}
C.~L. Epstein.
\newblock Envelopes of horospheres and weingarten surfaces in hyperbolic
  3-space.
\newblock 1984.

\bibitem{epstein-marden}
D.~B.~A. Epstein and A.~Marden.
\newblock Convex hulls in hyperbolic spaces, a theorem of {Sullivan}, and
  measured pleated surfaces.
\newblock In D.~B.~A. Epstein, editor, {\em Analytical and geometric aspects of
  hyperbolic space}, volume 111 of {\em L.M.S. Lecture Note Series}. Cambridge
  University Press, 1986.

\bibitem{FM2011}
B.~Farb and D.~Margalit.
\newblock {\em {Primer on Mapping Class Groups}}.
\newblock {Princeton Mathematical Series}, 2011.

\bibitem{gabai-meyerhoff-milley}
D.~Gabai, R.~Meyerhoff, and P.~Milley.
\newblock Minimum volume cusped hyperbolic three-manifolds.
\newblock {\em J. Amer. Math. Soc.}, 22(4):1157--1215, 2009.

\bibitem{Weeks}
D.~{Gabai}, R.~{Meyerhoff}, and P.~{Milley}.
\newblock {Minimum volume cusped hyperbolic three-manifolds}.
\newblock {\em Journal of the American Mathematical Society}, 22(4):1157--1215,
  Oct. 2009.

\bibitem{QTT}
N.~Gardiner, Frederick P.;~Lakic.
\newblock {\em Quasiconformal Teichm\"uller Theory}, volume~76 of {\em
  Mathematical Surveys and Monographs}.
\newblock American Mathematical Society, first edition, 2000.

\bibitem{graham-witten}
C.~R. Graham and E.~Witten.
\newblock Conformal anomaly of submanifold observables in {A}d{S}/{CFT}
  correspondence.
\newblock {\em Nuclear Phys. B}, 546(1-2):52--64, 1999.

\bibitem{He1976}
J.~Hempel.
\newblock {\em {3-Manifolds}}.
\newblock {Princeton University Press}, 1976.

\bibitem{Hubbard2016}
J.~H. Hubbard.
\newblock {\em Teichm{\"u}ller theory and applications to geometry, topology,
  and dynamics}, volume 1-2.
\newblock Matrix Editions, 2016.

\bibitem{klingenberg}
W.~P.~A. Klingenberg.
\newblock {\em Riemannian geometry}, volume~1 of {\em De Gruyter Studies in
  Mathematics}.
\newblock Walter de Gruyter \& Co., Berlin, second edition, 1995.

\bibitem{kojima-mcshane}
S.~Kojima and G.~McShane.
\newblock Normalized entropy versus volume for pseudo-{A}nosovs.
\newblock {\em Geom. Topol.}, 22(4):2403--2426, 2018.

\bibitem{Holography}
K.~Krasnov.
\newblock Holography and {R}iemann surfaces.
\newblock {\em Adv. Theor. Math. Phys.}, 4(4):929--979, 2000.

\bibitem{S2008}
K.~Krasnov and J.-M. Schlenker.
\newblock On the renormalized volume of hyperbolic 3-manifolds.
\newblock {\em Communications in Mathematical Physics}, pages 637--668, 2008.

\bibitem{kraus1932zusammenhang}
W.~Kraus.
\newblock {\"U}ber den zusammenhang einiger charakteristiken eines einfach
  zusammenh{\"a}ngenden bereiches mit der kreisabbildung.
\newblock {\em Mitteilungen des Mathematischen Seminars Giessen}, 21:1--28,
  1932.

\bibitem{Ma2016}
A.~Marden.
\newblock {\em {Hyperbolic Manifolds: An introduction in 2 and 3 dimensions}}.
\newblock {Cambridge University Press}, 2016.

\bibitem{Mar16}
B.~Martelli.
\newblock {An introduction to geometric topology}.
\newblock {\url{arXiv:1610.02592v1}}, 2016.

\bibitem{MT1998}
K.~Matsuzaki and M.~Taniguchi.
\newblock {\em {Hyperbolic Menifolds and Kleinian Groups}}.
\newblock {Oxford University Press}, 1998.

\bibitem{Ne1949}
Z.~Nehari.
\newblock The schwarzian derivative and schlicht functions.
\newblock {\em Bull. Amer. Math. Soc.}, 55:545--551, 1949.

\bibitem{scannell-wolf}
K.~P. Scannell and M.~Wolf.
\newblock The grafting map of {T}eichm\"uller space.
\newblock {\em J. Amer. Math. Soc.}, 15(4):893--927 (electronic), 2002.

\bibitem{compare}
J.-M. Schlenker.
\newblock The renormalized volume and the volume of the convex core of
  quasifuchsian manifolds.
\newblock {\em Math. Res. Lett.}, 20(4):773--786, 2013.
\newblock Corrected version available as arXiv:1109.6663v4.

\bibitem{volumes}
J.-M. Schlenker.
\newblock Volumes of quasifuchsian manifolds.
\newblock {\em Surveys in Differential Geometry}, 25(1):319--353, 2020.
\newblock arXiv:1903.09849.

\bibitem{averages}
J.-M. Schlenker and E.~Witten.
\newblock No ensemble averaging below the black hole threshold.
\newblock {\em J. High Energy Phys.}, pages Paper No. 143, 50, 2022.

\bibitem{TZ-schottky}
L.~Takhtajan and P.~Zograf.
\newblock On uniformization of {Riemann} surfaces and the {Weil-Petersson}
  metric on the {Teichm\"uller} and {Schottky} spaces.
\newblock {\em Mat. Sb.}, 132:303--320, 1987.
\newblock English translation in {\it Math. USSR Sb. 60:297-313, 1988}.

\bibitem{Takhtajan:2002cc}
L.~A. Takhtajan and L.-P. Teo.
\newblock Liouville action and {W}eil-{P}etersson metric on deformation spaces,
  global {K}leinian reciprocity and holography.
\newblock {\em Comm. Math. Phys.}, 239(1-2):183--240, 2003.

\bibitem{thurston-notes}
W.~P. Thurston.
\newblock Three-dimensional geometry and topology.
\newblock Originally notes of lectures at Princeton University, 1979. Recent
  version available on http://www.msri.org/publications/books/gt3m/, 1980.

\bibitem{VP2019}
F.~Vargas~Pallete.
\newblock Upper bounds on renormalized volume for schottky groups.
\newblock arxiv:1905.03303, 2019.

\bibitem{MD2019}
F.~Vargas~Pallete and D.~Mart{\'\i}nez-Granado.
\newblock Comparing hyperbolic and extremal lengths for shortest curves.
\newblock arxiv:1911.09078, 2019.

\end{thebibliography}
